\input amstex
\documentstyle{amsppt}
%
\catcode`@=11
\redefine\output@{%
  \def\break{\penalty-\@M}\let\par\endgraf
  \ifodd\pageno\global\hoffset=105pt\else\global\hoffset=8pt\fi  
  \shipout\vbox{%
    \ifplain@
      \let\makeheadline\relax \let\makefootline\relax
    \else
      \iffirstpage@ \global\firstpage@false
        \let\rightheadline\frheadline
        \let\leftheadline\flheadline
      \else
        \ifrunheads@ 
        \else \let\makeheadline\relax
        \fi
      \fi
    \fi
    \makeheadline \pagebody \makefootline}%
  \advancepageno \ifnum\outputpenalty>-\@MM\else\dosupereject\fi
}
\catcode`\@=\active
\nopagenumbers
\def\negskp{\hskip -2pt}
\def\tr{\operatorname{tr}}
\def\const{\operatorname{const}}
\def\Alpha{\operatorname{A}}
\def\vtrule{\vrule height 12pt depth 6pt}
\def\compos{\,\raise 1pt\hbox{$\sssize\circ$} \,}
\def\vtttrule{\vrule height 12pt depth 19pt}
\def\boxit#1#2{\vcenter{\hsize=122pt\offinterlineskip\hrule
  \line{\vtttrule\hss\vtop{\hsize=120pt\centerline{#1}\vskip 5pt
  \centerline{#2}}\hss\vtttrule}\hrule}}
\accentedsymbol\bbd{\kern 2pt\bar{\kern -2pt\bold d}}
\accentedsymbol\bulletH{\overset{\kern 2pt\sssize\bullet}\to H}
\accentedsymbol\circH{\overset{\kern 2pt\sssize\circ}\to H}
\accentedsymbol\bulletBH{\overset{\sssize\bullet}\to{\bold H}}
\accentedsymbol\circBH{\overset{\sssize\circ}\to{\bold H}}
\accentedsymbol\bulletd{\overset{\kern 3pt\sssize\bullet\kern -2pt}\to d}
\accentedsymbol\circd{\overset{\kern 3pt\sssize\circ\kern -2pt}\to d}
\accentedsymbol\bulbulgamma{\overset{\kern 1pt\sssize\bullet%
\bullet}\to\gamma}
\accentedsymbol\bulcircgamma{\overset{\kern 1pt\sssize\bullet%
\circ}\to\gamma}
\accentedsymbol\circbulgamma{\overset{\kern 1pt\sssize\circ%
\bullet}\to\gamma}
\accentedsymbol\circcircgamma{\overset{\kern 1pt\sssize\circ%
\circ}\to\gamma}
\def\blue#1{#1}
\catcode`#=11\def\diez{#}\catcode`#=6
\catcode`_=11\def\podcherkivanie{_}\catcode`_=8
\def\mycite#1{\cite{\blue{#1}}\immediate\special{ps:
     ShrHPSdict begin /ShrBORDERthickness 0 def}}

\def\mytag#1{%
    \tag#1}
\def\mythetag#1{\thetag{\blue{#1}}\immediate\special{ps:
     ShrHPSdict begin /ShrBORDERthickness 0 def}}
\def\myrefno#1{\no#1}
\def\myhref#1#2{\blue{#2}\immediate\special{ps:
     ShrHPSdict begin /ShrBORDERthickness 0 def}}
\def\myEarXivlink{\myhref{http://arXiv.org}{http:/\negskp/arXiv.org}}
\def\myGeoCities{\myhref{http://www.geocities.com}{GeoCities}}
\def\mytheorem#1{\csname proclaim\endcsname{Theorem #1}}
\def\mythetheorem#1{\blue{#1}\immediate\special{ps:
     ShrHPSdict begin /ShrBORDERthickness 0 def}}
\def\mylemma#1{\csname proclaim\endcsname{Lemma #1}}

\def\mycorollary#1{\csname proclaim\endcsname{Corollary #1}}

\def\mydefinition#1{\definition{Definition #1}}

\pagewidth{360pt}
\pageheight{606pt}
\topmatter
\title
Comparison of two formulas for metric connections 
in the bundle of Dirac spinors.
\endtitle
\author
R.~A.~Sharipov
\endauthor
\address 5 Rabochaya street, 450003 Ufa, Russia\newline
\vphantom{a}\kern 12pt Cell Phone: +7-(917)-476-93-48
\endaddress
\email \vtop to 30pt{\hsize=280pt\noindent
\myhref{mailto:r-sharipov\@mail.ru}
{r-sharipov\@mail.ru}\newline
\myhref{mailto:R\podcherkivanie Sharipov\@ic.bashedu.ru}
{R\_\hskip 1pt Sharipov\@ic.bashedu.ru}\vss}
\endemail
\urladdr
\vtop to 20pt{\hsize=280pt\noindent
\myhref{http://www.geocities.com/r-sharipov}
{http:/\negskp/www.geocities.com/r-sharipov}\newline
\myhref{http://www.freetextbooks.boom.ru/index.html}
{http:/\negskp/www.freetextbooks.boom.ru/index.html}\vss}
\endurladdr
\abstract
    Two explicit formulas for metric connections in the bundle 
of Dirac spinors are studied. Their equivalence is proved. The
explicit formula relating the spinor curvature tensor with the 
Riemann curvature tensor is rederived.
\endabstract
\endtopmatter
\loadbold
\loadeufb
\TagsOnRight
\document
\accentedsymbol\tbvartheta{\tilde{\overline{\boldsymbol\vartheta}}
\vphantom{\boldsymbol\vartheta}}

\rightheadtext{Comparison of two formulas for metric connections \dots}
\head
1. Basic notations and definitions. 
\endhead
    Dirac spinors play crucial role in modern particle physics. 
However, this crucial application of Dirac spinors is based mostly
on the special relativity, where the base manifold $M$ is the flat
Minkowski space. Passing to the general relativity, we get a 
little bit more complicated theory of spinors.\par
    Let $M$ be a {\it space-time} manifold of the general relativity. 
It is a four-dimensional orientable manifold equipped with a 
pseudo-Euclidean Minkowski-type metric $\bold g$ and with a 
{\it polarization}. The polarization of $M$ is responsible for 
distinguishing the {\it Future light cone\/} from the 
{\it Past light cone\/} at each point $p\in M$ (see \mycite{1} 
for more details). Let's denote by $DM$ the bundle of Dirac spinors 
over $M$ (see \mycite{2} and \mycite{3} for detailed description). 
In addition to the metric tensor $\bold g$ inherited from $M$, 
the Dirac bundle $DM$ is equipped with four other basic
spin-tensorial fields:
$$
\vcenter{\hsize 10cm
\offinterlineskip\settabs\+\indent
\vtrule
\hskip 1.2cm &\vtrule 
\hskip 5.2cm &\vtrule 
\hskip 2.8cm &\vtrule 
\cr\hrule 
\+\vtrule
\hfill\,Symbol\hfill&\vtrule
\hfill Name\hfill &\vtrule
\hfill Spin-tensorial\hfill &\vtrule\cr
\vskip -0.2cm
\+\vtrule
\hfill &\vtrule
\hfill \hfill&\vtrule
\hfill type\hfill&\vtrule\cr\hrule
\+\vtrule
\hfill $\bold g$\hfill&\vtrule
\hfill Metric tensor\hfill&\vtrule
\hfill $(0,0|0,0|0,2)$\hfill&\vtrule\cr\hrule
\+\vtrule
\hfill $\bold d$\hfill&\vtrule
\hfill Skew-symmetric metric tensor\hfill&\vtrule
\hfill $(0,2|0,0|0,0)$\hfill&\vtrule\cr\hrule
\+\vtrule
\hfill$\bold H$\hfill&\vtrule
\hfill Chirality operator\hfill&\vtrule
\hfill $(1,1|0,0|0,0)$\hfill&\vtrule\cr\hrule
\+\vtrule
\hfill$\bold D$\hfill&\vtrule
\hfill Dirac form\hfill&\vtrule
\hfill $(0,1|0,1|0,0)$\hfill&\vtrule\cr\hrule
\+\vtrule
\hfill$\boldsymbol\gamma$\hfill&\vtrule
\hfill Dirac $\gamma$-field\hfill&\vtrule
\hfill $(1,1|0,0|1,0)$\hfill&\vtrule\cr\hrule
}\quad
\mytag{1.1}
$$
The spin-tensorial type in the above table \mythetag{1.1} reflects
the number of indices in coordinate representation of fields. The
first two numbers are the numbers of upper and lower spinor indices, 
the second two numbers are the numbers of upper and lower conjugate
spinor indices, and the last two numbers are the numbers of upper 
and lower tensorial indices (they are also called spacial indices).
The metric tensor $\bold g$ is interpreted as a spin-tensorial field 
of the type $(0,0|0,0|0,2)$, i\.\,e\. it has no spinor indices and no
conjugate spinor indices, but has two lower spacial indices.\par
     The Dirac bundle is a complex bundle over a real manifold. 
For this reason spin-tensorial bundles produced from $DM$ are equipped 
with the involution of complex conjugation $\tau$ which exchanges 
spinor and conjugate spinor indices:
$$
\hskip -2em
\CD
@>\tau>>\\
\vspace{-4ex}
D^\alpha_\beta\bar D^\nu_\gamma T^m_n M@.
D^\nu_\gamma\bar D^\alpha_\beta T^m_n M.\\
\vspace{-4.2ex}
@<<\tau< 
\endCD
\mytag{1.2}
$$
Two fields $\bold g$ and $\bold D$ in \mythetag{1.1} are
real fields, i.\,e\. they are invariant with respect to
the involution of complex conjugation \mythetag{1.2}:
$$
\xalignat 2
&\tau(\bold g)=\bold g,
&&\tau(\bold D)=\bold D.
\endxalignat
$$
Other fields $\bold d$, $\bold H$, and $\boldsymbol\gamma$ in 
\mythetag{1.1} are not real fields. For them we denote
$$
\xalignat 3
&\bar{\boldsymbol\gamma}=\tau(\boldsymbol\gamma),
&&\bbd=\tau(\bold d),
&&\bar{\bold H}=\tau(\bold H).
\endxalignat 
$$
\mydefinition{1.1} A metric connection $(\Gamma,\Alpha,\bar{\Alpha})$ 
in $DM$ is a spinor connection which is real in the sense of the 
involution \mythetag{1.2} and concordant with $\bold d$ and 
$\boldsymbol\gamma$, i\.\,e\.
$$
\xalignat 3
&\tau\compos\nabla=\nabla\compos\tau,
&&\nabla\bold d=0,
&&\nabla\boldsymbol\gamma=0.
\ \qquad
\mytag{1.3}
\endxalignat
$$
\enddefinition
     Note that we use three symbols $(\Gamma,\Alpha,\bar{\Alpha})$ 
for denoting a spinor connection. It is because we use three different
types of connection components for three groups of indices when writing 
covariant derivatives in coordinates:
$$
\align
&\hskip -2em
\left.
\aligned
&\nabla_{\kern -1pt i}X^k=L_{\boldsymbol\Upsilon_i}X^k
+\shave{\sum^3_{j=0}}\Gamma^k_{ij}\,X^j,\\
&\nabla_{\kern -1pt i}X_k=L_{\boldsymbol\Upsilon_i}X_k
-\shave{\sum^3_{j=0}}\Gamma^j_{i\kern 1pt k}\,X_j,
\endaligned\right\}
\text{\ for spacial indices},
\mytag{1.4}\\
&\hskip -2em
\left.
\aligned
&\nabla_{\kern -1pt i}\psi^a=L_{\boldsymbol\Upsilon_i}\psi^a
+\shave{\sum^4_{b=1}}\Alpha^a_{i\kern 1pt b}\,\psi^b,\\
&\nabla_{\kern -1pt i}\psi_a=L_{\boldsymbol\Upsilon_i}\psi_a
-\shave{\sum^4_{b=1}}\Alpha^b_{i\kern 1pt a}\,\psi_b,
\endaligned\right\}
\text{\ for spinor indices},
\mytag{1.5}\\
&\hskip -2em
\left.
\aligned
&\nabla_{\kern -1pt i}\psi^{\bar a}=L_{\boldsymbol\Upsilon_i}
\psi^{\bar a}+\shave{\sum^4_{\bar b=1}}\bar{\Alpha}
\vphantom{\Alpha}^{\bar a}_{i\kern 1pt \bar b}\,\psi^{\bar b},\\
&\nabla_{\kern -1pt i}\psi_{\bar a}=L_{\boldsymbol\Upsilon_i}
\psi_{\bar a}-\shave{\sum^4_{\bar b=1}}\bar{\Alpha}
\vphantom{\Alpha}^{\bar b}_{i\kern 1pt \bar a}\,\psi_{\bar b},
\endaligned\right\}
\ \aligned
&\text{for conjugate}\\
&\text{spinor indices}.
\endaligned
\mytag{1.6}\\
\endalign
$$
In the case of a field of some mixed spin-tensorial type the 
formulas \mythetag{1.4}, \mythetag{1.5}, and \mythetag{1.6} are 
combined (see formula \thetag{7.10} in \mycite{3}).\par
     Note that in \mythetag{1.4}, \mythetag{1.5}, and \mythetag{1.6}
we have no partial derivatives. They are replaced by the derivatives
$L_{\boldsymbol\Upsilon_i}$ along four vector fields $\boldsymbol
\Upsilon_0,\,\boldsymbol\Upsilon_1,\,\boldsymbol\Upsilon_2,\,\boldsymbol
\Upsilon_3$ forming some local frame of the tangent bundle $TM$. The
spacial indices $i$, $j$, and $k$ in \mythetag{1.4}, \mythetag{1.5}, 
and \mythetag{1.6} are also relative to this frame. The spinor and
conjugate spinor indices $a$, $b$, $\bar a$, and $\bar b$ in these 
formulas are relative to some spinor frame $\boldsymbol\Psi_1,
\,\boldsymbol\Psi_2,\,\boldsymbol\Psi_3,\,\boldsymbol\Psi_4$ of 
the Dirac bundle $DM$. Local frames of the tangent bundle $TM$ are 
generalizations of local coordinates (they are also called non-holonomic
coordinates). Indeed, once some local coordinates $x^0,\,x^1,\,x^2,
\,x^3$ in $M$ are given, we have their associated frame of coordinate
vector fields:
$$
\xalignat 4
&\hskip -2em
\bold E_0=\frac{\partial}{\partial x^0},
&&\bold E_1=\frac{\partial}{\partial x^1},
&&\bold E_2=\frac{\partial}{\partial x^2},
&&\bold E_3=\frac{\partial}{\partial x^3}.
\quad
\mytag{1.7}
\endxalignat
$$
The local coordinates $x^0,\,x^1,\,x^2,\,x^3$ are called holonomic, 
since their associated vector fields \mythetag{1.7} commute with each 
other:
$$
\hskip -2em
[\bold E_i,\bold E_j]=0.
\mytag{1.8}
$$
Unlike \mythetag{1.8}, the vector fields $\boldsymbol\Upsilon_0,\,
\boldsymbol\Upsilon_1,\,\boldsymbol\Upsilon_2,\,\boldsymbol\Upsilon_3$
of a general (non-holonomic) frame do not commute. In this case we have
$$
\hskip -2em
[\boldsymbol\Upsilon_i,\boldsymbol\Upsilon_j]=\sum^3_{k=0}
c^{\,k}_{ij}\,\boldsymbol\Upsilon_k.
\mytag{1.9}
$$
As for spinor frames, they are always non-holonomic since they are
composed by spinor fields, while the commutator of spinor fields is 
not defined at all.\par
\mytheorem{1.1} Any metric connection $(\Gamma,\Alpha,\bar{\Alpha})$
is concordant with all of the basic spin tensorial fields $\bold g$, 
$\bold d$, $\bold H$, $\bold D$, and $\boldsymbol\gamma$ listed in the
table \mythetag{1.1}.
\endproclaim
\noindent The theorem~\mythetheorem{1.1} means that from \mythetag{1.3} it
follows that 
$$
\xalignat 3
&\nabla\bold g=0,
&&\nabla\bold H=0,
&&\nabla\bold D=0.
\ \qquad
\mytag{1.10}
\endxalignat
$$
Applying the first identity \mythetag{1.3} to \mythetag{1.10} and
to other identities \mythetag{1.3}, we derive
$$
\xalignat 3
&\nabla\bar{\boldsymbol\gamma}=0,
&&\nabla\bbd=0,
&&\nabla\bar{\bold H}=0.
\ \qquad
\endxalignat
$$
The general relativity (the Einstein's theory of gravity) is a theory
with zero torsion. Exactly for this case we have the following theorem.
\mytheorem{1.2} There is a unique metric connection $(\Gamma,\Alpha,
\bar{\Alpha})$ of the bundle of Dirac spinors $DM$ whose torsion is 
zero.
\endproclaim
    The metric connection with zero torsion is called the 
{\it Levi-Civita\/} connection. The proof of both 
theorems~\mythetheorem{1.1} and ~\mythetheorem{1.2} can be found 
in \mycite{3}. The first identity \mythetag{1.10} means that 
$\Gamma^k_{ij}$ in \mythetag{1.4} are the components of the
standard Levi-Civita connection. In a holonomic frame \mythetag{1.7}
they are given by the standard formula
$$
\hskip -2em
\Gamma^k_{ij}=\sum^3_{r=0}\frac{g^{\kern 0.5pt kr}}{2}
\left(\frac{\partial g_{rj}}{\partial x^i}+\frac{\partial g_{i
\kern 0.5ptr}}{\partial x^j}-\frac{\partial g_{ij}}{\partial x^r}
\right)\!.
\mytag{1.11}
$$
The quantities \mythetag{1.11} are symmetric, i\.\,e\. $\Gamma^k_{ij}
=\Gamma^k_{j\kern 1pt i}$. In the case of a non-holonomic frame
$\boldsymbol\Upsilon_0,\,\boldsymbol\Upsilon_1,\,\boldsymbol\Upsilon_2,
\,\boldsymbol\Upsilon_3$ the components of the Levi-Civita connection
are not symmetric: $\Gamma^k_{ij}-\Gamma^k_{j\kern 1pt i}=-c^{\,k}_{ij}$
with $c^{\,k}_{ij}$ taken from \mythetag{1.9}. They are given by the formula
$$
\hskip -2em
\gathered
\Gamma^k_{ij}=\sum^3_{r=0}\frac{g^{\kern 0.5pt kr}}{2}
\left(L_{\boldsymbol\Upsilon_{\!i}}\!(g_{rj})
+L_{\boldsymbol\Upsilon_{\!j}}\!(g_{i\kern 0.5pt r})
-L_{\boldsymbol\Upsilon_{\!r}}\!(g_{ij})\right)-\\
-\,\frac{c^{\,k}_{ij}}{2}
+\sum^3_{r=0}\sum^3_{s=0}\frac{c^{\,s}_{i\kern 0.5pt r}}{2}\,g^{kr}
\,g_{sj}+\sum^3_{r=0}\sum^3_{s=0}\frac{c^{\,s}_{j\kern 0.5ptr}}{2}
\,g^{kr}\,g_{s\kern 0.5pt i}.
\endgathered
\mytag{1.12}
$$\par
    Now let's proceed to the quantities $\Alpha^a_{i\kern 1pt b}$ and
$\bar{\Alpha}\vphantom{\Alpha}^{\bar a}_{i\kern 1pt \bar b}$ in the
formulas \mythetag{1.5} and \mythetag{1.6}. Due to the 
theorem~\mythetheorem{1.2} for a torsion-free connection they are
uniquely determined by the equalities \mythetag{1.3}. From the first
equality \mythetag{1.3} one easily derives that $\Alpha^a_{i\kern 1pt b}$
and $\bar{\Alpha}\vphantom{\Alpha}^{\bar a}_{i\kern 1pt \bar b}$ are
related to each other by virtue of the complex conjugation:
$$
\hskip -2em
\bar{\Alpha}\vphantom{\Alpha}^{\bar a}_{i\kern 1pt \bar b}
=\overline{\Alpha^{\bar a}_{i\kern 1pt \bar b}}.
\mytag{1.13}
$$
In one of my previous papers I have derived the following formula for 
the spinor components $\Alpha^a_{i\kern 1pt b}$ of the metric connection
(see formula \thetag{8.34} in \mycite{3}):
$$
\hskip -2em
\gathered
\Alpha^a_{ib}=\sum^4_{\alpha=1}\sum^4_{\beta=1}
\frac{L_{\boldsymbol\Upsilon_{\!i}}(\bulletd_{\kern 0.5pt \alpha\beta})
\,\bulletd^{\kern 0.5pt \beta\kern 0.5pt \alpha}\,\circH^a_b
+L_{\boldsymbol\Upsilon_{\!i}}(\circd_{\kern 0.5pt \alpha\beta})\,
\circd^{\kern 0.5pt \beta\kern 0.5pt \alpha}\,\bulletH^a_b}{4}\,+\\
+\sum^3_{m=0}\sum^3_{n=0}\sum^4_{\alpha=1}
\frac{L_{\boldsymbol\Upsilon_{\!i}}(\bulcircgamma^{\,\alpha}_{\!b
\kern 0.5pt m})\,g^{mn}\,\circbulgamma^{\,a}_{\!\alpha n}+
L_{\boldsymbol\Upsilon_{\!i}}(\circbulgamma^{\,\alpha}_{\!b
\kern 0.5pt m})\,g^{mn}\,\bulcircgamma^{\,a}_{\!\alpha n}}{4}\,-\\
-\sum^3_{m=0}\sum^3_{n=0}\sum^3_{r=0}\sum^4_{\alpha=1}
\frac{\Gamma^r_{im}\,\bulcircgamma^{\,\alpha}_{\!b\kern 0.5pt r}
\,g^{mn}\,\circbulgamma^{\,a}_{\!\alpha n}
+\Gamma^r_{im}\,\circbulgamma^{\,\alpha}_{\!b\kern 0.5pt r}
\,g^{mn}\,\bulcircgamma^{\,a}_{\!\alpha n}}{4}.
\endgathered\qquad
\mytag{1.14}
$$
However, in some other papers there are much more simple formulas for 
spinor connections. I choose the formula \thetag{5} from \mycite{4} 
for comparing it with \mythetag{1.14}. Being transformed to our 
notations, this formula looks like
$$
\hskip -2em
\aligned
A^a_{ib}&=\frac{1}{4}\sum^3_{m=0}\sum^3_{n=0}\sum^3_{r=0}
\sum^4_{\alpha=1}\Gamma^r_{im}\,\gamma^\alpha_{bn}\,g^{mn}
\,\gamma^a_{\alpha r}.
\endaligned
\mytag{1.15}
$$    
Our main goal is to compare the formulas \mythetag{1.14} and
\mythetag{1.15}. Then we calculate the curvature tensors associated
with these formulas.
\head
2. Comparison of the formulas (1.14) and (1.15).
\endhead
    The formula \mythetag{1.15} has no differentiations at all. This
means that it is written for special frames where the components of
the basic fields $\bold d$ and $\boldsymbol\gamma$ are constants.
For this reason we could omit the terms with derivatives 
$L_{\boldsymbol\Upsilon_{\!i}}$ in \mythetag{1.14} for the comparison
purposes. However, we shall not do it. We shall transform the
formula \mythetag{1.14} to a form similar to \mythetag{1.15}
preserving all of its terms. As a result in each step we shall have
a formula applicable for arbitrary frames.\par
     Another feature of the formula \mythetag{1.14}, in contrast to
\mythetag{1.15}, is that it uses special notations with circles and
bullets. These notations were introduced in \mycite{3}. Now we need 
to reproduce them. Remember that the chirality operator $\bold H$ 
is a diagonalizable operator with two eigenvalues $\lambda_1=\lambda_2
=1$ and two other eigenvalues $\lambda_3=\lambda_4=-1$. Therefore, 
$\bold H^2=1$. Due to this equality one can define two projection 
operators 
$$
\xalignat 2
&\hskip -2em
\bulletBH=\frac{1+\bold H}{2},
&&\circBH=\frac{1-\bold H}{2}
\mytag{2.1}
\endxalignat
$$
(see \thetag{8.1} in \mycite{3}). The components $\circH^a_b$ and
$\bulletH^a_b$ of these two projection operators \mythetag{2.1} are
used in the formula \mythetag{1.14}. Other special symbols in
\mythetag{1.14} are defined as 
$$
\xalignat 2
&\hskip -2em
\bulletd^{\kern 0.5pt \beta\kern 0.5pt \alpha}
=\sum^4_{a=1}\sum^4_{b=1}
\bulletH^\beta_b\ d^{\kern 0.5pt b\kern 0.5pt a}
\,\bulletH^\alpha_a,
&&\circd^{\kern 0.5pt \beta\kern 0.5pt \alpha}
=\sum^4_{a=1}\sum^4_{b=1}
\circH^\beta_b\ d^{\kern 0.5pt b\kern 0.5pt a}
\,\circH^\alpha_a,\qquad\\
\vspace{-1.5ex}
&&&\mytag{2.2}\\
\vspace{-1.5ex}
&\hskip -2em
\bulletd_{\kern 0.5pt \alpha\beta}=\sum^4_{a=1}\sum^4_{b=1}
\bulletH^a_\alpha\ d_{\kern 0.5pt ab}\,\bulletH^b_\beta,
&&\circd_{\kern 0.5pt \alpha\beta}=\sum^4_{a=1}\sum^4_{b=1}
\circH^a_\alpha\ d_{\kern 0.5pt ab}\,\circH^b_\beta\qquad
\endxalignat
$$
(see \thetag{8.19} in \mycite{3}). Here $d_{\kern 0.5pt ab}$
and $d^{\kern 0.5pt b\kern 0.5pt a}$ are the components of two
mutually inverse skew-symmetric matrices. The first of them
represents the spinor metric tensor $\bold d$ (see table 
\mythetag{1.1}) and the second one corresponds to its dual
metric tensor which is denoted by the same symbol $\bold d$.
Similarly, we have
$$
\xalignat 2
&\hskip -2em
\circbulgamma^{\,\alpha}_{\!\beta m}=\sum^4_{a=1}\sum^4_{b=1}
\circH^\alpha_a\,\bulletH^b_\beta\,\gamma^{\,a}_{b\kern 0.5pt m},
&&\bulcircgamma^{\,\alpha}_{\!\beta m}=\sum^4_{a=1}\sum^4_{b=1}
\bulletH^\alpha_a\,\circH^b_\beta\,\gamma^{\,a}_{b\kern 0.5pt m}
\quad
\mytag{2.3}
\endxalignat
$$
(see \thetag{8.5} in \mycite{3}). The projection operators 
\mythetag{2.1} obey some commutation and anticommutation relationships
with $\bold d$ and $\boldsymbol\gamma$:
$$
\xalignat 2
&\hskip -2em
\sum^4_{b=1}d_{\alpha b}\,\bulletH^b_\beta
=\sum^4_{a=1}\bulletH^a_\alpha\,d_{a\beta},
&&\sum^4_{a=1}d^{\kern 0.5pt\beta\kern 0.5pt a}\,\bulletH^\alpha_a
=\sum^4_{b=1}\bulletH^\beta_b\,d^{\kern 0.5pt b\kern 0.5pt \alpha},\\
\vspace{-1.5ex}
&&&\mytag{2.4}\\
\vspace{-1.5ex}
&\hskip -2em
\sum^4_{b=1}d_{\alpha b}\,\circH^b_\beta
=\sum^4_{a=1}\circH^a_\alpha\,d_{a\beta},
&&\sum^4_{a=1}d^{\kern 0.5pt\beta\kern 0.5pt a}\,\circH^\alpha_a
=\sum^4_{b=1}\circH^\beta_b\,d^{\kern 0.5pt b\kern 0.5pt \alpha},\\
&\hskip -2em
\sum^4_{b=1}\gamma^{\,\alpha}_{b\kern 0.5pt m}\,\bulletH^b_\beta
=\sum^4_{a=1}\circH^\alpha_a\,\gamma^{\,a}_{\beta\kern 0.5pt m},
&&\sum^4_{b=1}\gamma^{\,\alpha}_{b\kern 0.5pt m}\,\circH^b_\beta
=\sum^4_{a=1}\bulletH^\alpha_a\,\gamma^{\,a}_{\beta\kern 0.5pt m}
\quad
\mytag{2.5}
\endxalignat
$$
(see \thetag{6.16} and \thetag{6.17} in \mycite{3}). Due to the
relationships \mythetag{2.4} and \mythetag{2.5} the formulas 
\mythetag{2.2} and \mythetag{2.3} are written as follows:
$$
\allowdisplaybreaks
\align
&\hskip -2em
\aligned
\bulletd_{\kern 0.5pt \alpha\beta}=\sum^4_{b=1}d_{\alpha b}
\,\bulletH^b_\beta &=\sum^4_{a=1}\bulletH^a_\alpha\,d_{a\beta},\\
\bulletd^{\kern 0.5pt \beta\kern 0.5pt \alpha}
&=\sum^4_{a=1}d^{\kern 0.5pt\beta\kern 0.5pt a}\,\bulletH^\alpha_a
=\sum^4_{b=1}\bulletH^\beta_b\,d^{\kern 0.5pt b\kern 0.5pt \alpha},
\endaligned
\mytag{2.6}\\
&\hskip -2em
\aligned
\circd_{\kern 0.5pt \alpha\beta}=\sum^4_{b=1}d_{\alpha b}
\,\circH^b_\beta &=\sum^4_{a=1}\circH^a_\alpha\,d_{a\beta},\\
\circd^{\kern 0.5pt \beta\kern 0.5pt \alpha}
&=\sum^4_{a=1}d^{\kern 0.5pt\beta\kern 0.5pt a}\,\circH^\alpha_a
=\sum^4_{b=1}\circH^\beta_b\,d^{\kern 0.5pt b\kern 0.5pt \alpha},
\endaligned
\mytag{2.7}\\
&\hskip -2em
\aligned
\circbulgamma^{\,\alpha}_{\!\beta m}
=\sum^4_{b=1}\gamma^{\,\alpha}_{b\kern 0.5pt m}\,\bulletH^b_\beta
&=\sum^4_{a=1}\circH^\alpha_a\,\gamma^{\,a}_{\beta\kern 0.5pt m},\\
\bulcircgamma^{\,\alpha}_{\!\beta m}
&=\sum^4_{b=1}\gamma^{\,\alpha}_{b\kern 0.5pt m}\,\circH^b_\beta
=\sum^4_{a=1}\bulletH^\alpha_a\,\gamma^{\,a}_{\beta\kern 0.5pt m}.
\endaligned
\mytag{2.8}
\endalign
$$
Now we apply \mythetag{2.6}, \mythetag{2.7}, and \mythetag{2.8}
for to transform the formula \mythetag{1.14}. Using the formulas 
\mythetag{2.8}, we derive the following relationship:
$$
\hskip -2em
\gathered
\sum^4_{\alpha=1}
\Bigl(\bulcircgamma^{\,\alpha}_{\!b\kern 0.5pt r}
\,\circbulgamma^{\,a}_{\!\alpha n}+
\circbulgamma^{\,\alpha}_{\!b\kern 0.5pt r}
\,\bulcircgamma^{\,a}_{\!\alpha n}\Bigr)=
\sum^4_{\alpha=1}\Bigl(\circbulgamma^{\,a}_{\!\alpha n}\,
\bulcircgamma^{\,\alpha}_{\!b\kern 0.5pt r}
+\bulcircgamma^{\,a}_{\!\alpha n}
\,\circbulgamma^{\,\alpha}_{\!b\kern 0.5pt r}\Bigr)=\\
=\sum^4_{\alpha=1}\sum^4_{c=1}\sum^4_{d=1}
\Bigl(\gamma^{\,a}_{c\kern 0.5pt n}\,\bulletH^c_\alpha
\,\bulletH^\alpha_d\,\gamma^{\,d}_{b\kern 0.5pt r}
+\gamma^{\,a}_{c\kern 0.5pt n}\,\circH^c_\alpha
\,\circH^\alpha_d\,\gamma^{\,d}_{b\kern 0.5pt r}\Bigr).
\endgathered
\mytag{2.9}
$$
Remember that $\bulletBH$ and $\circBH$ are projection operators
complementary to each other, i\.\,e\. $\bulletBH^2=\bulletBH$,
$\circBH^2=\circBH$, and $\bulletBH+\circBH=1$. Therefore, we have 
$$
\hskip -2em
\sum^4_{\alpha=1}\Bigl(\bulletH^c_\alpha\,\bulletH^\alpha_d
+\circH^c_\alpha\,\circH^\alpha_d\Bigr)=\bulletH^c_d+\circH^c_d
=\delta^c_d.
\mytag{2.10}
$$
Applying this formula \mythetag{2.10} to \mythetag{2.9} we
derive
$$
\hskip -2em
\sum^4_{\alpha=1}
\Bigl(\bulcircgamma^{\,\alpha}_{\!b\kern 0.5pt r}
\,\circbulgamma^{\,a}_{\!\alpha n}+
\circbulgamma^{\,\alpha}_{\!b\kern 0.5pt r}
\,\bulcircgamma^{\,a}_{\!\alpha n}\Bigr)
=\sum^4_{c=1}\gamma^{\,a}_{c\kern 0.5pt n}
\,\gamma^{\,c}_{b\kern 0.5pt r}.
\mytag{2.11}
$$
Using the formula \mythetag{2.11}, we can transform the last term
of \mythetag{1.14} as follows:
$$
\hskip -2em
\gathered
\sum^3_{m=0}\sum^3_{n=0}\sum^3_{r=0}\sum^4_{\alpha=1}
\frac{\Gamma^r_{im}\,\bulcircgamma^{\,\alpha}_{\!b\kern 0.5pt r}
\,g^{mn}\,\circbulgamma^{\,a}_{\!\alpha n}
+\Gamma^r_{im}\,\circbulgamma^{\,\alpha}_{\!b\kern 0.5pt r}
\,g^{mn}\,\bulcircgamma^{\,a}_{\!\alpha n}}{4}=\\
\vspace{2ex}
=\sum^3_{m=0}\sum^3_{n=0}\sum^3_{r=0}\sum^4_{c=1}
\frac{\Gamma^r_{im}\,\gamma^{\,c}_{b\kern 0.5pt r}
\,g^{mn}\,\gamma^{\,a}_{cn}}{4}.
\endgathered
\mytag{2.12}
$$\par
     Now let's proceed to the second term in the right hand side
of the formula \mythetag{1.14}. Applying the formulas \mythetag{2.8},
we derive
$$
\gather
\sum^3_{m=0}\sum^3_{n=0}\sum^4_{\alpha=1}
\Bigl(\,L_{\boldsymbol\Upsilon_i}(\bulcircgamma^{\,\alpha}_{\!b
\kern 0.5pt m})\,g^{mn}\,\circbulgamma^{\,a}_{\!\alpha n}
+L_{\boldsymbol\Upsilon_i}(\circbulgamma^{\,\alpha}_{\!b
\kern 0.5pt m})\,g^{mn}\,\bulcircgamma^{\,a}_{\!\alpha n}\Bigr)=\\
=\sum^3_{m=0}\sum^3_{n=0}\sum^4_{\alpha=1}\sum^4_{c=1}\sum^4_{d=1}
\gamma^{\,a}_{c\kern 0.5pt n}\left(\bulletH^c_\alpha
\,L_{\boldsymbol\Upsilon_i}(\bulletH^\alpha_d\,\gamma^{\,d}_{b
\kern 0.5pt m})+\circH^c_\alpha\,L_{\boldsymbol\Upsilon_i}
(\circH^\alpha_d\,\gamma^{\,d}_{b\kern 0.5pt m})\right)g^{mn}=\\
=\sum^3_{m=0}\sum^3_{n=0}\sum^4_{\alpha=1}\sum^4_{c=1}\sum^4_{d=1}
\gamma^{\,a}_{c\kern 0.5pt n}\left(\bulletH^c_\alpha
\,L_{\boldsymbol\Upsilon_i}(\bulletH^\alpha_d)
+\circH^c_\alpha\,L_{\boldsymbol\Upsilon_i}(\circH^\alpha_d)\right)
g^{mn}\,\gamma^{\,d}_{b\kern 0.5pt m}\,+\\
+\sum^3_{m=0}\sum^3_{n=0}\sum^4_{\alpha=1}\sum^4_{c=1}\sum^4_{d=1}
\gamma^{\,a}_{c\kern 0.5pt n}\left(\bulletH^c_\alpha
\,\bulletH^\alpha_d+\circH^c_\alpha\,\circH^\alpha_d\right)
L_{\boldsymbol\Upsilon_i}(\gamma^{\,d}_{b\kern 0.5pt m})
\,g^{mn}.
\endgather
$$
In order to continue our calculations we need the identity
\thetag{6.21} from \mycite{3}:
$$
\hskip -2em
\gathered
\sum^3_{m=0}\sum^3_{n=0}\gamma^{\,d}_{b\kern 0.5pt m}\ g^{mn}
\,\gamma^{\,a}_{c\kern 0.5pt n}
=\delta^d_c\ \delta^a_b
-H^d_c\ H^a_b\,+\\
+\,d^{\kern 0.5pt da}\,d_{\kern 0.5pt bc}-\sum^4_{r=1}\sum^4_{s=1}
H^d_r\,d^{\kern 0.5pt ra}\,d_{\kern 0.5pt bs}\,H^s_c.
\endgathered
\mytag{2.13}
$$
Applying \mythetag{2.10} and \mythetag{2.13} to the above formula,
we derive
$$
\gather
\sum^3_{m=0}\sum^3_{n=0}\sum^4_{\alpha=1}
\Bigl(\,L_{\boldsymbol\Upsilon_i}(\bulcircgamma^{\,\alpha}_{\!b
\kern 0.5pt m})\,g^{mn}\,\circbulgamma^{\,a}_{\!\alpha n}
+L_{\boldsymbol\Upsilon_i}(\circbulgamma^{\,\alpha}_{\!b
\kern 0.5pt m})\,g^{mn}\,\bulcircgamma^{\,a}_{\!\alpha n}\Bigr)=\\
=\sum^3_{m=0}\sum^3_{n=0}\sum^4_{\alpha=1}
L_{\boldsymbol\Upsilon_i}(\gamma^{\,\alpha}_{b\kern 0.5pt m})
\,g^{mn}\,\gamma^{\,a}_{\alpha\kern 0.5pt n}\,+\\
=\sum^4_{\alpha=1}\sum^4_{c=1}\sum^4_{d=1}\bulletH^c_\alpha
\,L_{\boldsymbol\Upsilon_i}(\bulletH^\alpha_d)
\left(\delta^d_c\ \delta^a_b-H^d_c\ H^a_b\,+\vphantom{\sum^4_{r=1}}
\qquad\qquad\right.\\
\left.\qquad\qquad\qquad\qquad+\,d^{\kern 0.5pt da}\,d_{\kern 0.5pt bc}
-\sum^4_{r=1}\sum^4_{s=1}
H^d_r\,d^{\kern 0.5pt ra}\,d_{\kern 0.5pt bs}\,H^s_c\right)+\\
+\sum^4_{\alpha=1}\sum^4_{c=1}\sum^4_{d=1}\circH^c_\alpha
\,L_{\boldsymbol\Upsilon_i}(\circH^\alpha_d)
\left(\delta^d_c\ \delta^a_b-H^d_c\ H^a_b\,+\vphantom{\sum^4_{r=1}}
\qquad\qquad\right.\\
\left.\qquad\qquad\qquad\qquad+\,d^{\kern 0.5pt da}\,d_{\kern 0.5pt
bc}-\sum^4_{r=1}\sum^4_{s=1}
H^d_r\,d^{\kern 0.5pt ra}\,d_{\kern 0.5pt bs}\,H^s_c\right).
\endgather
$$
Now remember that $\bulletBH\compos\bold H=\bold H\compos\bulletBH
=\bulletBH$ and $\circBH\compos\bold H=\bold H\compos\circBH
=-\circBH$. These formulas are easily derived from \mythetag{2.1}.
Applying them, we find
$$
\gather
\sum^3_{m=0}\sum^3_{n=0}\sum^4_{\alpha=1}
\Bigl(\,L_{\boldsymbol\Upsilon_i}(\bulcircgamma^{\,\alpha}_{\!b
\kern 0.5pt m})\,g^{mn}\,\circbulgamma^{\,a}_{\!\alpha n}
+L_{\boldsymbol\Upsilon_i}(\circbulgamma^{\,\alpha}_{\!b
\kern 0.5pt m})\,g^{mn}\,\bulcircgamma^{\,a}_{\!\alpha n}\Bigr)=\\
=\sum^3_{m=0}\sum^3_{n=0}\sum^4_{\alpha=1}
L_{\boldsymbol\Upsilon_i}(\gamma^{\,\alpha}_{b\kern 0.5pt m})
\,g^{mn}\,\gamma^{\,a}_{\alpha\kern 0.5pt n}
+\sum^4_{\alpha=1}\sum^4_{d=1}\bulletH^d_\alpha
\,L_{\boldsymbol\Upsilon_i}(\bulletH^\alpha_d)
\left(\delta^a_b-H^a_b\right)+\\
+\sum^4_{\alpha=1}\sum^4_{c=1}\sum^4_{d=1}d_{\kern 0.5pt bc}
\,\bulletH^c_\alpha\,L_{\boldsymbol\Upsilon_i}(\bulletH^\alpha_d)
\sum^4_{r=1}\left(\delta^d_r-H^d_r\right)d^{\kern 0.5pt ra}
+\sum^4_{\alpha=1}\sum^4_{d=1}\circH^d_\alpha
\,L_{\boldsymbol\Upsilon_i}(\circH^\alpha_d)\,\times\\
\times\,\left(\delta^a_b+H^a_b\right)
+\sum^4_{\alpha=1}\sum^4_{c=1}\sum^4_{d=1}d_{\kern 0.5pt bc}
\,\circH^c_\alpha\,L_{\boldsymbol\Upsilon_i}(\circH^\alpha_d)
\sum^4_{r=1}\left(\delta^d_r+H^d_r\right)d^{\kern 0.5pt ra}.
\endgather
$$
Note that $\delta^a_b-H^a_b=2\,\circH^a_b$ and $\delta^a_b+H^a_b
=2\,\bulletH^a_b$. Similarly, $\delta^d_r-H^d_r=2\,\circH^d_r$
and $\delta^d_r+H^d_r=2\,\bulletH^d_r$. Moreover, we apply the 
following obvious formula
$$
L_{\boldsymbol\Upsilon_i}(\gamma^{\,\alpha}_{b\kern 0.5pt m})
\,g^{mn}=L_{\boldsymbol\Upsilon_i}(\gamma^{\,\alpha}_{b\kern 0.5pt m}
\,g^{mn})-\gamma^{\,\alpha}_{b\kern 0.5pt m}\,L_{\boldsymbol\Upsilon_i}
(g^{mn}).
$$

Therefore, the above 
formula is transformed to the following one:
$$
\gather
\sum^3_{m=0}\sum^3_{n=0}\sum^4_{\alpha=1}
\Bigl(\,L_{\boldsymbol\Upsilon_i}(\bulcircgamma^{\,\alpha}_{\!b
\kern 0.5pt m})\,g^{mn}\,\circbulgamma^{\,a}_{\!\alpha n}
+L_{\boldsymbol\Upsilon_i}(\circbulgamma^{\,\alpha}_{\!b
\kern 0.5pt m})\,g^{mn}\,\bulcircgamma^{\,a}_{\!\alpha n}\Bigr)=\\
=\sum^3_{m=0}\sum^3_{n=0}\sum^4_{\alpha=1}
L_{\boldsymbol\Upsilon_i}(\gamma^{\,\alpha}_{b\kern 0.5pt m}\,g^{mn})
\,\gamma^{\,a}_{\alpha\kern 0.5pt n}
-\sum^3_{m=0}\sum^3_{n=0}\sum^4_{\alpha=1}
\gamma^{\,\alpha}_{b\kern 0.5pt m}\,L_{\boldsymbol\Upsilon_i}(g^{mn})
\,\gamma^{\,a}_{\alpha\kern 0.5pt n}\,+\\
+\,2\sum^4_{\alpha=1}\sum^4_{d=1}\bulletH^d_\alpha
\,L_{\boldsymbol\Upsilon_i}(\bulletH^\alpha_d)
\,\circH^a_b+2\sum^4_{\alpha=1}\sum^4_{c=1}\sum^4_{d=1}
\sum^4_{r=1}d_{\kern 0.5pt bc}\,\bulletH^c_\alpha
\,L_{\boldsymbol\Upsilon_i}(\bulletH^\alpha_d)
\,\circH^d_r\,d^{\kern 0.5pt ra}\,+\\
+\,2\sum^4_{\alpha=1}\sum^4_{d=1}\circH^d_\alpha
\,L_{\boldsymbol\Upsilon_i}(\circH^\alpha_d)\,\bulletH^a_b
+2\sum^4_{\alpha=1}\sum^4_{c=1}\sum^4_{d=1}\sum^4_{r=1}
d_{\kern 0.5pt bc}\,\circH^c_\alpha
\,L_{\boldsymbol\Upsilon_i}(\circH^\alpha_d)
\,\bulletH^d_r\,d^{\kern 0.5pt ra}.
\endgather
$$
For the derivative $L_{\boldsymbol\Upsilon_i}(g^{mn})$ in the
above formula we write
$$
\hskip -2em
L_{\boldsymbol\Upsilon_i}(g^{mn})=
-\sum^3_{s=0}\Gamma^m_{is}\,g^{sn}
-\sum^3_{s=0}\Gamma^n_{is}\,g^{ms}.
\mytag{2.14}
$$
This formula is easily derived from $\nabla_{\kern -1pt i}g^{mn}=0$.
Applying \mythetag{2.14}, we get
$$
\gathered
\sum^3_{m=0}\sum^3_{n=0}\sum^4_{\alpha=1}
\Bigl(\,L_{\boldsymbol\Upsilon_i}(\bulcircgamma^{\,\alpha}_{\!b
\kern 0.5pt m})\,g^{mn}\,\circbulgamma^{\,a}_{\!\alpha n}
+L_{\boldsymbol\Upsilon_i}(\circbulgamma^{\,\alpha}_{\!b
\kern 0.5pt m})\,g^{mn}\,\bulcircgamma^{\,a}_{\!\alpha n}\Bigr)=\\
=\sum^3_{m=0}\sum^3_{n=0}\sum^4_{\alpha=1}
L_{\boldsymbol\Upsilon_i}(\gamma^{\,\alpha}_{b\kern 0.5pt m}
\,g^{mn})\,\gamma^{\,a}_{\alpha\kern 0.5pt n}
+\sum^3_{m=0}\sum^3_{n=0}\sum^4_{\alpha=1}\sum^3_{s=0}
\gamma^{\,\alpha}_{b\kern 0.5pt m}\,\Gamma^m_{is}\,g^{sn}
\,\gamma^{\,a}_{\alpha\kern 0.5pt n}\,+\\
+\sum^3_{m=0}\sum^3_{n=0}\sum^4_{\alpha=1}\sum^3_{s=0}
\gamma^{\,\alpha}_{b\kern 0.5pt m}\,\Gamma^n_{is}\,g^{ms}
\,\gamma^{\,a}_{\alpha\kern 0.5pt n}
+2\sum^4_{\alpha=1}\sum^4_{d=1}\bulletH^d_\alpha
\,L_{\boldsymbol\Upsilon_i}(\bulletH^\alpha_d)
\,\circH^a_b\,+\\
+\,2\sum^4_{\alpha=1}\sum^4_{c=1}\sum^4_{d=1}
\sum^4_{r=1}d_{\kern 0.5pt bc}\,\bulletH^c_\alpha
\,L_{\boldsymbol\Upsilon_i}(\bulletH^\alpha_d)
\,\circH^d_r\,d^{\kern 0.5pt ra}
+2\sum^4_{\alpha=1}\sum^4_{d=1}\circH^d_\alpha\times\\
\times\,L_{\boldsymbol\Upsilon_i}(\circH^\alpha_d)
\,\bulletH^a_b+2\sum^4_{\alpha=1}\sum^4_{c=1}\sum^4_{d=1}
\sum^4_{r=1}d_{\kern 0.5pt bc}\,\circH^c_\alpha
\,L_{\boldsymbol\Upsilon_i}(\circH^\alpha_d)
\,\bulletH^d_r\,d^{\kern 0.5pt ra}.
\endgathered\qquad
\mytag{2.15}
$$\par
     The last step in our calculations is to transform the first term
in the right hand side of \mythetag{1.14}. Applying \mythetag{2.6} and
\mythetag{2.7} to it, we get
$$
\gather
\sum^4_{\alpha=1}\sum^4_{\beta=1}\Bigl(
L_{\boldsymbol\Upsilon_{\!i}}(\bulletd_{\kern 0.5pt \alpha\beta})
\,\bulletd^{\kern 0.5pt \beta\kern 0.5pt \alpha}\,\circH^a_b
+L_{\boldsymbol\Upsilon_{\!i}}(\circd_{\kern 0.5pt \alpha\beta})\,
\circd^{\kern 0.5pt \beta\kern 0.5pt \alpha}\,\bulletH^a_b\Bigr)=\\
=\sum^4_{\alpha=1}\sum^4_{\beta=1}\sum^4_{c=1}\sum^4_{d=1}
L_{\boldsymbol\Upsilon_{\!i}}(\bulletH^c_\alpha\,d_{c\beta})
\,d^{\kern 0.5pt \beta\kern 0.5pt d}\,\bulletH^\alpha_d\,\circH^a_b
+\sum^4_{\alpha=1}\sum^4_{\beta=1}\sum^4_{c=1}\sum^4_{d=1}
L_{\boldsymbol\Upsilon_{\!i}}(\circH^c_\alpha\,d_{c\beta})\times\\
\times\,d^{\kern 0.5pt \beta\kern 0.5pt d}\,\circH^\alpha_d
\,\bulletH^a_b=
\sum^4_{\alpha=1}\sum^4_{d=1}
L_{\boldsymbol\Upsilon_{\!i}}(\bulletH^d_\alpha)
\,\bulletH^\alpha_d\,\circH^a_b
+\sum^4_{\alpha=1}\sum^4_{\beta=1}\sum^4_{d=1}
L_{\boldsymbol\Upsilon_{\!i}}(d_{\alpha\beta})
\,d^{\kern 0.5pt \beta\kern 0.5pt d}\,\bulletH^\alpha_d\,\times\\
\times\,\circH^a_b
+\sum^4_{\alpha=1}\sum^4_{d=1}
L_{\boldsymbol\Upsilon_{\!i}}(\circH^d_\alpha)
\,\circH^\alpha_d\,\bulletH^a_b
+\sum^4_{\alpha=1}\sum^4_{\beta=1}\sum^4_{d=1}
L_{\boldsymbol\Upsilon_{\!i}}(d_{\alpha\beta})
\,d^{\kern 0.5pt \beta\kern 0.5pt d}\,\circH^\alpha_d
\,\bulletH^a_b.
\endgather
$$
Some terms in the above sum are zero. Indeed, we have
$$
\gathered
\sum^4_{\alpha=1}\sum^4_{d=1}
L_{\boldsymbol\Upsilon_{\!i}}(\bulletH^d_\alpha)
\,\bulletH^\alpha_d=\frac{1}{4}\sum^4_{\alpha=1}\sum^4_{d=1}
L_{\boldsymbol\Upsilon_{\!i}}(H^d_\alpha)
\left(\delta^\alpha_d+H^\alpha_d\right)=\\
=\frac{1}{4}\tr(L_{\boldsymbol\Upsilon_{\!i}}(\bold H))
+\frac{1}{4}\tr(L_{\boldsymbol\Upsilon_{\!i}}(\bold H)
\compos\bold H)
=\frac{1}{4}L_{\boldsymbol\Upsilon_{\!i}}(\tr\bold H)+
\frac{1}{8}L_{\boldsymbol\Upsilon_{\!i}}(\tr\bold H^2).
\endgathered
$$
But we know that $\tr\bold H=0$ and $\bold H^2=1$, which means
$\tr\bold H^2=4$. Therefore, we have the following two relationships:
$$
\hskip -2em
\aligned
&\sum^4_{\alpha=1}\sum^4_{d=1}
L_{\boldsymbol\Upsilon_{\!i}}(\bulletH^d_\alpha)
\,\bulletH^\alpha_d=0,\\
&\sum^4_{\alpha=1}\sum^4_{d=1}
L_{\boldsymbol\Upsilon_{\!i}}(\circH^d_\alpha)
\,\circH^\alpha_d=0.
\endaligned
\mytag{2.16}
$$
The second relationship \mythetag{2.16} is derived in a quite similar 
way as the first one. Applying these relationships, we continue our
previous calculations and obtain 
$$
\gathered
\sum^4_{\alpha=1}\sum^4_{\beta=1}
L_{\boldsymbol\Upsilon_{\!i}}(\bulletd_{\kern 0.5pt \alpha\beta})
\,\bulletd^{\kern 0.5pt \beta\kern 0.5pt \alpha}\,\circH^a_b
+L_{\boldsymbol\Upsilon_{\!i}}(\circd_{\kern 0.5pt \alpha\beta})\,
\circd^{\kern 0.5pt \beta\kern 0.5pt \alpha}\,\bulletH^a_b=\\
=\sum^4_{\alpha=1}\sum^4_{\beta=1}\sum^4_{d=1}\Bigl(
L_{\boldsymbol\Upsilon_{\!i}}(d_{\alpha\beta})
\,d^{\kern 0.5pt \beta\kern 0.5pt d}\,\bulletH^\alpha_d\,
\circH^a_b+L_{\boldsymbol\Upsilon_{\!i}}(d_{\alpha\beta})
\,d^{\kern 0.5pt \beta\kern 0.5pt d}\,\circH^\alpha_d
\,\bulletH^a_b\Bigr).
\endgathered
\mytag{2.17}
$$
Applying \mythetag{2.16} to \mythetag{2.15} we can simplify it either:
$$
\gathered
\sum^3_{m=0}\sum^3_{n=0}\sum^4_{\alpha=1}
\Bigl(\,L_{\boldsymbol\Upsilon_i}(\bulcircgamma^{\,\alpha}_{\!b
\kern 0.5pt m})\,g^{mn}\,\circbulgamma^{\,a}_{\!\alpha n}
+L_{\boldsymbol\Upsilon_i}(\circbulgamma^{\,\alpha}_{\!b
\kern 0.5pt m})\,g^{mn}\,\bulcircgamma^{\,a}_{\!\alpha n}\Bigr)=\\
=\sum^3_{m=0}\sum^3_{n=0}\sum^4_{\alpha=1}
L_{\boldsymbol\Upsilon_i}(\gamma^{\,\alpha}_{b\kern 0.5pt m}
\,g^{mn})\,\gamma^{\,a}_{\alpha\kern 0.5pt n}
+\sum^3_{m=0}\sum^3_{n=0}\sum^4_{\alpha=1}\sum^3_{s=0}
\gamma^{\,\alpha}_{b\kern 0.5pt m}\,\Gamma^m_{is}\,g^{sn}
\,\gamma^{\,a}_{\alpha\kern 0.5pt n}\,+\\
+\sum^3_{m=0}\sum^3_{n=0}\sum^4_{\alpha=1}\sum^3_{s=0}
\gamma^{\,\alpha}_{b\kern 0.5pt m}\,\Gamma^n_{is}\,g^{ms}
\,\gamma^{\,a}_{\alpha\kern 0.5pt n}
+2\sum^4_{\alpha=1}\sum^4_{c=1}\sum^4_{d=1}
\sum^4_{r=1}d_{\kern 0.5pt bc}\,\bulletH^c_\alpha\,\times\\
\times\,L_{\boldsymbol\Upsilon_i}(\bulletH^\alpha_d)
\,\circH^d_r\,d^{\kern 0.5pt ra}
+2\sum^4_{\alpha=1}\sum^4_{c=1}\sum^4_{d=1}
\sum^4_{r=1}d_{\kern 0.5pt bc}\,\circH^c_\alpha
\,L_{\boldsymbol\Upsilon_i}(\circH^\alpha_d)
\,\bulletH^d_r\,d^{\kern 0.5pt ra}.
\endgathered\qquad
\mytag{2.18}
$$\par
Moreover, note that $\bulletBH$ and $\circBH$ are projectors, 
i\.\,e\. $\bulletBH\compos\bulletBH=\bulletBH$ and $\circBH\compos
\circBH=\circBH$. Differentiating these formulas, we derive the 
following  identities:
$$
\aligned
&\sum^4_{\alpha=1}\sum^4_{d=1}\bulletH^c_\alpha
\,L_{\boldsymbol\Upsilon_i}(\bulletH^\alpha_d)
\,\circH^d_r=\sum^4_{\alpha=1}\bulletH^c_\alpha
\,L_{\boldsymbol\Upsilon_i}(\bulletH^\alpha_r)
=\sum^4_{d=1}\,L_{\boldsymbol\Upsilon_i}(\bulletH^c_d)
\,\circH^d_r,\\
&\sum^4_{\alpha=1}\sum^4_{d=1}\circH^c_\alpha
\,L_{\boldsymbol\Upsilon_i}(\circH^\alpha_d)
\,\bulletH^d_r=\sum^4_{\alpha=1}\circH^c_\alpha
\,L_{\boldsymbol\Upsilon_i}(\circH^\alpha_r)
=\sum^4_{d=1}\,L_{\boldsymbol\Upsilon_i}(\circH^c_d)
\,\bulletH^d_r.
\endaligned
\quad
\mytag{2.19}
$$
These formulas will be used below. Now we substitute \mythetag{2.17},
\mythetag{2.18}, and \mythetag{2.12} into the formula \mythetag{1.14}.
Meanwhile we apply \mythetag{2.19} to \mythetag{2.18}. As a result 
we get
$$
\gathered
\Alpha^a_{ib}
=\sum^4_{\alpha=1}\sum^4_{\beta=1}\sum^4_{d=1}\frac{
L_{\boldsymbol\Upsilon_{\!i}}(d_{\alpha\beta})
\,d^{\kern 0.5pt \beta\kern 0.5pt d}\,\bulletH^\alpha_d\,
\circH^a_b+L_{\boldsymbol\Upsilon_{\!i}}(d_{\alpha\beta})
\,d^{\kern 0.5pt \beta\kern 0.5pt d}\,\circH^\alpha_d
\,\bulletH^a_b}{4}\,+\\
+\sum^4_{c=1}\sum^4_{d=1}
\sum^4_{r=1}\frac{d_{\kern 0.5pt bc}
\,L_{\boldsymbol\Upsilon_i}(\bulletH^c_d)
\,\circH^d_r\,d^{\kern 0.5pt ra}}{2}
+\sum^4_{c=1}\sum^4_{d=1}
\sum^4_{r=1}\frac{d_{\kern 0.5pt bc}
\,L_{\boldsymbol\Upsilon_i}(\circH^c_d)
\,\bulletH^d_r\,d^{\kern 0.5pt ra}}{2}+\\
+\sum^3_{m=0}\sum^3_{n=0}\sum^4_{\alpha=1}
\frac{L_{\boldsymbol\Upsilon_i}(\gamma^{\,\alpha}_{b\kern 0.5pt m}
\,g^{mn})\,\gamma^{\,a}_{\alpha\kern 0.5pt n}}{4}
+\sum^3_{m=0}\sum^3_{n=0}\sum^4_{\alpha=1}\sum^3_{s=0}
\frac{\gamma^{\,\alpha}_{b\kern 0.5pt m}\,\Gamma^n_{is}\,g^{ms}
\,\gamma^{\,a}_{\alpha\kern 0.5pt n}}{4}.
\endgathered\qquad
\mytag{2.20}
$$
Note that the last term in \mythetag{2.20} coincides with 
\mythetag{1.15}, other terms contain derivatives 
$L_{\boldsymbol\Upsilon_i}$. Due to this observation we
can formulate the following result.
\mytheorem{2.1} The formulas \mythetag{1.14} and \mythetag{1.15}
represent the same metric connection for Dirac spinors. The formula
\mythetag{1.15} applies to special frame where the components of the 
basic fields listed in the table \mythetag{1.1} are constants. The
formula \mythetag{1.14} is a general formula applicable to all frames.
It can be written as \mythetag{2.20}.
\endproclaim
\head
3. Further transformations of the formula (1.14).
\endhead
     The third line in the formula \mythetag{2.20} has no symbols with
circle and bullet.  However, previous two lines still have such symbols.
We use the formulas \mythetag{2.1} to remove bullets and circles from
the formula \mythetag{2.20} at all:
$$
\gather
\hskip -2em
\gathered
\bulletH^\alpha_d\,\circH^a_b+\circH^\alpha_d\,\bulletH^a_b
=\frac{1}{2}\,(\delta^\alpha_d+H^\alpha_d)\,\circH^a_b
+\frac{1}{2}\,(\delta^\alpha_d-H^\alpha_d)\,\bulletH^a_b=\\
=\frac{1}{2}\,\delta^\alpha_d\,(\circH^a_b+\bulletH^a_b)
+\frac{1}{2}\,H^\alpha_d\,(\circH^a_b-\bulletH^a_b)
=\frac{1}{2}\,\delta^\alpha_d\,\delta^a_b
-\frac{1}{2}\,H^\alpha_d\,H^a_b,
\endgathered
\mytag{3.1}\\
\vspace{2ex}
\hskip -2em
\gathered
L_{\boldsymbol\Upsilon_i}(\bulletH^c_d)\,\circH^d_r
+L_{\boldsymbol\Upsilon_i}(\circH^c_d)\,\bulletH^d_r
=\frac{1}{2}\,L_{\boldsymbol\Upsilon_i}(\delta^c_d+H^c_d)
\,\circH^d_r\,+\\
+\,\frac{1}{2}\,L_{\boldsymbol\Upsilon_i}(\delta^c_d-H^c_d)
\,\bulletH^d_r=
\frac{1}{2}\,L_{\boldsymbol\Upsilon_i}(H^c_d)\,(\circH^d_r
-\bulletH^d_r)=-\frac{1}{2}\,L_{\boldsymbol\Upsilon_i}(H^c_d)
\,H^d_r.
\endgathered\quad
\mytag{3.2}
\endgather
$$    
Applying \mythetag{3.1} and \mythetag{3.2} to \mythetag{2.20},
we obtain
$$
\gathered
\Alpha^a_{ib}
=\sum^4_{\alpha=1}\sum^4_{\beta=1}\frac{
L_{\boldsymbol\Upsilon_{\!i}}(d_{\alpha\beta})
\,d^{\kern 0.5pt \beta\kern 0.5pt\alpha}}{8}\,\delta^a_b
-\sum^4_{\alpha=1}\sum^4_{\beta=1}\sum^4_{d=1}\frac{
L_{\boldsymbol\Upsilon_{\!i}}(d_{\alpha\beta})
\,d^{\kern 0.5pt \beta\kern 0.5pt d}\,H^\alpha_d}{8}\,H^a_b\,-\\
-\sum^4_{c=1}\sum^4_{d=1}
\sum^4_{r=1}\frac{d_{\kern 0.5pt bc}
\,L_{\boldsymbol\Upsilon_i}(H^c_d)
\,H^d_r\,d^{\kern 0.5pt ra}}{4}\,+
\sum^3_{m=0}\sum^3_{n=0}\sum^4_{\alpha=1}
\frac{L_{\boldsymbol\Upsilon_i}(\gamma^{\,\alpha}_{b\kern 0.5pt m}
\,g^{mn})}{4}\,\times\\
\times\,\gamma^{\,a}_{\alpha\kern 0.5pt n}
+\sum^3_{m=0}\sum^3_{n=0}\sum^4_{\alpha=1}\sum^3_{s=0}
\frac{\gamma^{\,\alpha}_{b\kern 0.5pt m}\,\Gamma^n_{is}\,g^{ms}
\,\gamma^{\,a}_{\alpha\kern 0.5pt n}}{4}.
\endgathered\qquad
\mytag{3.3}
$$\par
    The theorem~\mythetheorem{1.1} says that the condition 
$\nabla\bold H=0$ follows from the conditions \mythetag{1.3}.
Therefore, one can expect that the derivatives \pagebreak
$L_{\boldsymbol\Upsilon_i}(H^c_d)$ are expressed through the
derivatives of the components of $\bold d$ and
$\boldsymbol\gamma$. It is really so and the calculations 
of $L_{\boldsymbol\Upsilon_i}(H^c_d)$ are quite similar to
those in proving the theorem~7.3 in \mycite{3}. These calculations
are based on the formula \mythetag{2.13}. For the sake of brevity
we omit them and give the ultimate result only. Here is the formula 
for $L_{\boldsymbol\Upsilon_i}(H^c_d)$:
$$
\gathered
L_{\boldsymbol\Upsilon_i}(H^c_d)=\sum^4_{\alpha=1}
\sum^4_{r=1}\frac{H^c_r\,d^{\kern 0.5pt r\alpha}
\,L_{\boldsymbol\Upsilon_i}(d_{\kern 0.5pt\alpha d})}{6}
-\sum^4_{\beta=1}\sum^4_{s=1}\frac{d^{\kern 0.5pt c\beta}
\,L_{\boldsymbol\Upsilon_i}(d_{\kern 0.5pt\beta s})
\,H^s_d}{6}+\\
+\sum^4_{\alpha=1}\sum^4_{\beta=1}
\frac{L_{\boldsymbol\Upsilon_i}(d^{\kern 0.5pt c\alpha}
\,d_{\kern 0.5pt \beta d})\,H^\beta_\alpha}{6}
-\sum^4_{\alpha=1}\sum^4_{\beta=1}\sum^3_{m=0}\sum^3_{n=0}
\frac{L_{\boldsymbol\Upsilon_i}(\gamma^{\,c}_{\beta\kern 0.5pt n}
\ g^{mn}\,\gamma^{\,\alpha}_{d\kern 0.5pt m})\,H^\beta_\alpha}{6}.
\endgathered\quad
\mytag{3.4}
$$
It is clear that substituting \mythetag{3.4} into \mythetag{3.3} will
make this formula more huge and complicated. Therefore, we stop our
transformations of the formula \mythetag{1.14} at this point assuming 
that \mythetag{3.3} is the most simple formula for the spinor components 
of the torsion-free metric connection $(\Gamma,\Alpha,\bar{\Alpha})$.
\par
\head
4. Special frames and curvature spin-tensors.
\endhead
    There are four types of special frames in the bundle of Dirac 
spinors $DM$. They are considered in \mycite{3}. A frame 
$\boldsymbol\Psi_1,\,\boldsymbol\Psi_2,\,\boldsymbol\Psi_3,
\,\boldsymbol\Psi_4$ of any one of these four types in $DM$ is 
canonically associated with some definite frame $\boldsymbol\Upsilon_0,
\,\boldsymbol\Upsilon_1,\,\boldsymbol\Upsilon_2,
\,\boldsymbol\Upsilon_3$ in the tangent bundle TM. The frame types 
association is given by the following diagram:
$$
\hskip -2em
\aligned
&\boxit{Canonically orthonormal}{chiral frames}\to
\boxit{Positively polarized}{right orthonormal frames}\\
&\boxit{$P$-reverse}{anti-chiral frames}\to
\boxit{Positively polarized}{left orthonormal frames}\\
&\boxit{$T$-reverse}{anti-chiral frames}\to
\boxit{Negatively polarized}{right orthonormal frames}\\
&\boxit{$PT$-reverse}{chiral frames}\to
\boxit{Negatively polarized}{left orthonormal frames}
\endaligned
\mytag{4.1}
$$
For the sake of certainty we choose some {\it canonically orthonormal
chiral frame\/} $\boldsymbol\Psi_1,\,\boldsymbol\Psi_2,\,\boldsymbol
\Psi_3,\,\boldsymbol\Psi_4$ in $DM$. According to the diagram 
\mythetag{4.1}, it is associated with some {\it positively polarized
right orthonormal frame\/} $\boldsymbol\Upsilon_0,\,\boldsymbol
\Upsilon_1,\,\boldsymbol\Upsilon_2,\,\boldsymbol\Upsilon_3$ in $TM$.
Then $g_{ij}=g(\boldsymbol\Upsilon_i,\boldsymbol\Upsilon_j)$ and for
the components of both metric tensors we have
$$
\hskip -2em
g_{ij}=g^{ij}=\Vmatrix\format \l&\quad\r&\quad\r&\quad\r\\
1 &0 &0 &0\\ 0 &-1 &0 &0\\ 0 &0 &-1 &0\\ 0 &0 &0 &-1
\endVmatrix.
\mytag{4.2}
$$
The formula \mythetag{4.2} means that $\boldsymbol\Upsilon_0,\,
\boldsymbol\Upsilon_1,\,\boldsymbol\Upsilon_2,\,\boldsymbol\Upsilon_3$
is an orthonormal frame in $TM$. Moreover, it is a right frame regarding
to the orientation in $M$. It is positively polarized, i\.\,e\. 
$\boldsymbol\Upsilon_0$ is a time-like unit vector directed to the
future.\par
     Canonically orthonormal chiral frames are simultaneously 
orthonormal, chiral, and self-adjoint frames in $DM$. The 
orthonormality of our frame $\boldsymbol\Psi_1,\,\boldsymbol\Psi_2,
\,\boldsymbol\Psi_3,\,\boldsymbol\Psi_4$ means that the components
of the metric tensor $\bold d$ are given by the following matrix:
$$
\hskip -2em
d_{ij}=
d(\boldsymbol\Psi_i,\boldsymbol\Psi_j)
=\Vmatrix 0 & 1 & 0 & 0\\-1 & 0 & 0 & 0\\
0 & 0 & 0 & -1\\0 & 0 & 1 & 0\endVmatrix.
\mytag{4.3}
$$\par
    Chiral frames in the bundle of Dirac spinors $DM$ are those 
for which the chirality operator $\bold H$ is given by the matrix
$$
\hskip -2em
H^{\kern 0.5pti}_{\kern -0.5pt j}=\Vmatrix 1 & 0 & 0 & 0\\0 & 1 & 0 & 0\\
0 & 0 & -1 & 0\\0 & 0 & 0 & -1\endVmatrix.
\mytag{4.4}
$$\par
     And finally, self-adjoint frames in the bundle of Dirac spinors 
$DM$ are those for which the Dirac form $\bold D$ is given by the
matrix
$$
\hskip -2em
D_{i\bar j}=D(\boldsymbol\Psi_{\bar j},\boldsymbol\Psi_i)
=\Vmatrix 0 & 0 & 1 & 0\\0 & 0 & 0 & 1\\
1 & 0 & 0 & 0\\0 & 1 & 0 & 0\endVmatrix.
\mytag{4.5}
$$\par
    Our choice is a canonically orthonormal chiral frame $\boldsymbol
\Psi_1,\,\boldsymbol\Psi_2,\,\boldsymbol\Psi_3,\,\boldsymbol\Psi_4$ in
$DM$ and its associated positively polarized right orthonormal frame 
$\boldsymbol\Upsilon_0,\,\boldsymbol \Upsilon_1,\,\boldsymbol
\Upsilon_2,\,\boldsymbol\Upsilon_3$ in $TM$. Therefore, in our case 
the conditions \mythetag{4.2}, \mythetag{4.3}, \mythetag{4.4}, and
\mythetag{4.5} are fulfilled simultaneously. The components of 
Dirac's $\gamma$-field are uniquely fixed by our choice of frames.
They are usually collected into four matrices:
$$
\xalignat 2
&\hskip -2em
\gamma_0=\Vmatrix 0&0&1&0\\0&0&0&1\\1&0&0&0\\0&1&0&0\endVmatrix,
&&\gamma_1=\Vmatrix 0&0&0&1\\0&0&1&0\\0&-1&0&0\\-1&0&0&0\endVmatrix,
\quad\\
\vspace{-1.5ex}
&&&\mytag{4.6}\\
\vspace{-1.5ex}
&\hskip -2em
\gamma_2=\Vmatrix 0&0&0&-i\\0&0&i&0\\0&i&0&0\\-i&0&0&0\endVmatrix,
&&\gamma_3=\Vmatrix 0&0&1&0\\0&0&0&-1\\-1&0&0&0\\0&1&0&0\endVmatrix.
\quad
\endxalignat
$$
The matrices \mythetag{4.6} are enumerated by the spacial index $k$
of $\gamma^{\,a}_{b\kern 0.5pt k}$. Other two indices $a$ and $b$ 
represent the position of an element within the matrix $\gamma_k$,
the index $a$ being the row number and the index $b$ being the column
number.\par
     Looking at the formulas \mythetag{4.2}, \mythetag{4.3}, 
\mythetag{4.4}, \mythetag{4.5}, and \mythetag{4.6}, we see that the
components of all basic fields $\bold g$, $\bold d$, $\bold H$,
$\bold D$, and $\boldsymbol\gamma$ are constants. It means that our
choice of frames is that very case where the formula \mythetag{1.15}
is applicable and where the formula \mythetag{1.14} reduces to
\mythetag{1.15}. This choice of frames is convenient for calculating
the curvature tensors. The first of them is the Riemann curvature 
tensor $\bold R$. Its components are given by the formula:
$$
R^p_{qij}=L_{\boldsymbol\Upsilon_i}(\Gamma^p_{\!j\,q})
-L_{\boldsymbol\Upsilon_j}(\Gamma^p_{\!i\,q})
+\sum^3_{h=0}\left(\Gamma^p_{\!i\,h}\,\Gamma^h_{\!j\,q}
-\Gamma^p_{\!j\,h}\,\Gamma^h_{\!i\,q}\right)
-\sum^3_{k=0}c^{\,k}_{ij}\,\Gamma^p_{kq}
\qquad
\mytag{4.7}
$$
(see \thetag{6.27} in \mycite{5}). Here $\Gamma^k_{ij}$ are the spacial
components of the metric connection $(\Gamma,\Alpha,\bar{\Alpha})$.
They are given by the formula \mythetag{1.12}. Apart from 
\mythetag{4.7}, there are two other curvature tensors $\eufb R$ and
$\bar{\eufb R}$. Their components are given by the formulas
$$
\align
&\hskip -4em
\goth R^p_{qij}=L_{\boldsymbol\Upsilon_i}(\Alpha^p_{j\,q})
-L_{\boldsymbol\Upsilon_j}(\Alpha^p_{i\,q})
+\sum^4_{h=1}\left(\Alpha^p_{i\,h}\,\Alpha^{\!h}_{j\,q}
-\Alpha^p_{j\,h}\,\Alpha^{\!h}_{i\,q}\right)
-\sum^3_{k=0}c^{\,k}_{ij}\,\Alpha^p_{kq},
\mytag{4.8}\\
&\hskip -4em
\bar{\goth R}^p_{qij}=L_{\boldsymbol\Upsilon_i}
(\bar{\Alpha}\vphantom{\Alpha}^p_{j\,q})
-L_{\boldsymbol\Upsilon_j}
(\bar{\Alpha}\vphantom{\Alpha}^p_{i\,q})
+\sum^4_{h=1}\left(\bar{\Alpha}\vphantom{\Alpha}^p_{i\,h}
\,\bar{\Alpha}\vphantom{\Alpha}^{\!h}_{j\,q}-\bar{\Alpha}
\vphantom{\Alpha}^p_{j\,h}\,\bar{\Alpha}
\vphantom{\Alpha}^{\!h}_{i\,q}\right)
-\sum^3_{k=0}c^{\,k}_{ij}\,\bar{\Alpha}
\vphantom{\Alpha}^p_{kq}
\mytag{4.9}
\endalign
$$
(compare with \thetag{6.25} and \thetag{6.26} in \mycite{5}).
Applying \mythetag{1.13} to \mythetag{4.8} and \mythetag{4.9}
and taking into account that $\Gamma^k_{ij}$ and $c^{\,k}_{ij}$
are purely real functions, we get
$$
\hskip -2em
\bar{\goth R}^p_{qij}=\overline{\goth R^p_{qij}}.
\mytag{4.10}
$$
In a coordinate-free form the relationship \mythetag{4.10} is
written as 
$$
\hskip -2em
\bar{\eufb R}=\tau(\eufb R),
\mytag{4.11}
$$
while the Riemann curvature tensor $\bold R$ is a purely real field,
i\.\,e\.
$$
\hskip -2em
\tau(\bold R)=\bold R.
\mytag{4.12}
$$\par
     Keeping in mind that we deal with the special frames
$\boldsymbol\Psi_1,\,\boldsymbol\Psi_2,\,\boldsymbol\Psi_3,
\,\boldsymbol\Psi_4$ and $\boldsymbol\Upsilon_0,\,\boldsymbol
\Upsilon_1,\,\boldsymbol\Upsilon_2,\,\boldsymbol\Upsilon_3$, 
let's apply the formula \mythetag{1.15} to \mythetag{4.8}.
As a result we get
$$
\gathered
\goth R^p_{qij}=\sum^3_{m=0}\sum^3_{n=0}\sum^3_{r=0}
\sum^4_{\alpha=1}\frac{L_{\boldsymbol\Upsilon_i}(\Gamma^r_{jm})
-L_{\boldsymbol\Upsilon_j}(\Gamma^r_{im})}{4}
\,\gamma^\alpha_{qn}\,g^{mn}\,\gamma^p_{\alpha r}+\\
\hskip -3em +\sum^4_{h=1}\sum^3_{m=0}\sum^3_{n=0}\sum^3_{r=0}
\sum^4_{\alpha=1}\frac{\Alpha^p_{i\,h}\,\Gamma^r_{jm}
\,\gamma^\alpha_{qn}\,g^{mn}\,\gamma^h_{\alpha r}}{4}\,-\\
\hskip 3em -\sum^4_{h=1}\sum^3_{m=0}\sum^3_{n=0}\sum^3_{r=0}
\sum^4_{\alpha=1}\frac{\Alpha^h_{i\,q}\,\Gamma^r_{jm}
\,\gamma^\alpha_{hn}\,g^{mn}\,\gamma^p_{\alpha r}}{4}\,-\\
\hskip -4em -\sum^3_{k=0}\sum^3_{m=0}\sum^3_{n=0}\sum^3_{r=0}
\sum^4_{\alpha=1}\frac{c^{\,k}_{ij}}{4}\,\Gamma^r_{km}
\,\gamma^\alpha_{qn}\,g^{mn}\,\gamma^p_{\alpha r}.
\endgathered\qquad
\mytag{4.13}
$$
Since $\gamma^p_{\alpha r}=\const$, the concordance condition 
$\nabla\boldsymbol\gamma=0$ from \mythetag{1.3} is written as
$$
\nabla_{\!i}\,\gamma^p_{\alpha r}=\sum^4_{h=1}\Alpha^p_{i\,h}
\,\gamma^h_{\alpha r}-\sum^4_{h=1}\Alpha^h_{i\,\alpha}
\,\gamma^p_{h\kern 0.5pt r}-\sum^3_{s=0}\Gamma^s_{i\kern 0.5pt r}
\,\gamma^p_{\alpha s}=0.
$$
From this identity we derive 
$$
\hskip -2em
\sum^4_{h=1}\Alpha^p_{i\,h}\,\gamma^h_{\alpha r}
=\sum^4_{h=1}\Alpha^h_{i\,\alpha}\,\gamma^p_{h\kern 0.5pt r}
+\sum^3_{s=0}\Gamma^s_{i\kern 0.5pt r}\,\gamma^p_{\alpha s}.
\mytag{4.14}
$$
Applying \mythetag{4.14} to the second term in the right hand 
side of \mythetag{4.13}, we write it as 
$$
\gathered
\goth R^p_{qij}=\sum^3_{m=0}\sum^3_{n=0}\sum^3_{r=0}
\sum^4_{\alpha=1}\frac{L_{\boldsymbol\Upsilon_i}(\Gamma^r_{jm})
-L_{\boldsymbol\Upsilon_j}(\Gamma^r_{im})}{4}
\,\gamma^\alpha_{qn}\,g^{mn}\,\gamma^p_{\alpha r}+\\
\hskip -3em +\sum^4_{h=1}\sum^3_{m=0}\sum^3_{n=0}\sum^3_{r=0}
\sum^4_{\alpha=1}\frac{\Alpha^h_{i\,\alpha}\,\Gamma^r_{jm}
\,\gamma^\alpha_{qn}\,g^{mn}\,\gamma^p_{h\kern 0.5pt r}}{4}\,+\\
\hskip -3em +\sum^3_{s=0}\sum^3_{m=0}\sum^3_{n=0}\sum^3_{r=0}
\sum^4_{\alpha=1}\frac{\Gamma^s_{i\kern 0.5pt r}\,\Gamma^r_{jm}
\,\gamma^\alpha_{qn}\,g^{mn}\,\gamma^p_{\alpha s}}{4}\,+\\
\hskip 3em -\sum^4_{h=1}\sum^3_{m=0}\sum^3_{n=0}\sum^3_{r=0}
\sum^4_{\alpha=1}\frac{\Alpha^h_{i\,q}\,\Gamma^r_{jm}
\,\gamma^\alpha_{hn}\,g^{mn}\,\gamma^p_{\alpha r}}{4}\,-\\
\hskip -4em -\sum^3_{k=0}\sum^3_{m=0}\sum^3_{n=0}\sum^3_{r=0}
\sum^4_{\alpha=1}\frac{c^{\,k}_{ij}}{4}\,\Gamma^r_{km}
\,\gamma^\alpha_{qn}\,g^{mn}\,\gamma^p_{\alpha r}.
\endgathered\qquad
\mytag{4.15}
$$
Now we use the following formula equivalent to \mythetag{4.14}:
$$
\hskip -2em
\sum^4_{\alpha=1}\Alpha^h_{i\,\alpha}\,\gamma^\alpha_{qn}
=\sum^4_{\alpha=1}\Alpha^\alpha_{i\,q}\,\gamma^h_{\alpha\kern 0.5pt n}
+\sum^3_{s=0}\Gamma^s_{i\kern 0.5pt n}\,\gamma^h_{qs}.
\mytag{4.16}
$$
Applying \mythetag{4.16} to the second term in the right hand
side of \mythetag{4.15}, we write it as
$$
\gathered
\goth R^p_{qij}=\sum^3_{m=0}\sum^3_{n=0}\sum^3_{r=0}
\sum^4_{\alpha=1}\frac{L_{\boldsymbol\Upsilon_i}(\Gamma^r_{jm})
-L_{\boldsymbol\Upsilon_j}(\Gamma^r_{im})}{4}
\,\gamma^\alpha_{qn}\,g^{mn}\,\gamma^p_{\alpha r}+\\
\hskip -3em +\sum^4_{h=1}\sum^3_{m=0}\sum^3_{n=0}\sum^3_{r=0}
\sum^4_{\alpha=1}\frac{\Alpha^\alpha_{i\,q}\,\Gamma^r_{jm}
\,\gamma^h_{\alpha\kern 0.5pt n}\,g^{mn}\,
\gamma^p_{h\kern 0.5pt r}}{4}\,+\\
\hskip 3em +\sum^4_{h=1}\sum^3_{m=0}\sum^3_{n=0}\sum^3_{r=0}
\sum^4_{\alpha=1}\frac{\Gamma^s_{i\kern 0.5pt n}\,\Gamma^r_{jm}
\,\gamma^h_{qs}\,g^{mn}\,\gamma^p_{h\kern 0.5pt r}}{4}\,+\\
\hskip -3em +\sum^3_{s=0}\sum^3_{m=0}\sum^3_{n=0}\sum^3_{r=0}
\sum^4_{\alpha=1}\frac{\Gamma^s_{i\kern 0.5pt r}\,\Gamma^r_{jm}
\,\gamma^\alpha_{qn}\,g^{mn}\,\gamma^p_{\alpha s}}{4}\,+\\
\hskip 3em -\sum^4_{h=1}\sum^3_{m=0}\sum^3_{n=0}\sum^3_{r=0}
\sum^4_{\alpha=1}\frac{\Alpha^h_{i\,q}\,\Gamma^r_{jm}
\,\gamma^\alpha_{hn}\,g^{mn}\,\gamma^p_{\alpha r}}{4}\,-\\
\hskip -4em -\sum^3_{k=0}\sum^3_{m=0}\sum^3_{n=0}\sum^3_{r=0}
\sum^4_{\alpha=1}\frac{c^{\,k}_{ij}}{4}\,\Gamma^r_{km}
\,\gamma^\alpha_{qn}\,g^{mn}\,\gamma^p_{\alpha r}.
\endgathered\qquad
\mytag{4.17}
$$
By means of the formal exchanges of summation indices 
$h\leftrightarrow\alpha$ we find that the second and the fifth 
terms in the right hand side of \mythetag{4.17} do cancel each
other. In order to transform the third term in the right hand 
side of \mythetag{4.17} we use the concordance condition
$\nabla\bold g=0$ from \mythetag{1.10}. Since $g^{ms}=\const$,
this condition yields
$$
\nabla_{\!i}g^{ms}=\sum^3_{n=0}\Gamma^m_{i\kern 0.5pt n}\,
g^{ns}+\sum^3_{n=0}\Gamma^s_{i\kern 0.5pt n}\,g^{mn}=0.
$$
This formula can be rewritten in the following way:
$$
\hskip -2em
\sum^3_{n=0}\Gamma^s_{i\kern 0.5pt n}\,g^{mn}=
-\sum^3_{n=0}\Gamma^m_{i\kern 0.5pt n}\,g^{ns}.
\mytag{4.18}
$$
Applying the formula \mythetag{4.18} to the third term in the
right hand side of \mythetag{4.17} and canceling the second
and the fifth terms there, we get
$$
\gathered
\goth R^p_{qij}=\sum^3_{m=0}\sum^3_{n=0}\sum^3_{r=0}
\sum^4_{\alpha=1}\frac{L_{\boldsymbol\Upsilon_i}(\Gamma^r_{jm})
-L_{\boldsymbol\Upsilon_j}(\Gamma^r_{im})}{4}
\,\gamma^\alpha_{qn}\,g^{mn}\,\gamma^p_{\alpha r}-\\
\hskip 3em -\sum^4_{h=1}\sum^3_{m=0}\sum^3_{n=0}\sum^3_{r=0}
\sum^4_{\alpha=1}\frac{\Gamma^m_{i\kern 0.5pt n}\,\Gamma^r_{jm}
\,\gamma^h_{qs}\,g^{ns}\,\gamma^p_{h\kern 0.5pt r}}{4}\,+\\
\hskip -3em +\sum^3_{s=0}\sum^3_{m=0}\sum^3_{n=0}\sum^3_{r=0}
\sum^4_{\alpha=1}\frac{\Gamma^s_{i\kern 0.5pt r}\,\Gamma^r_{jm}
\,\gamma^\alpha_{qn}\,g^{mn}\,\gamma^p_{\alpha s}}{4}\,-\\
\hskip 3em -\sum^3_{k=0}\sum^3_{m=0}\sum^3_{n=0}\sum^3_{r=0}
\sum^4_{\alpha=1}\frac{c^{\,k}_{ij}}{4}\,\Gamma^r_{km}
\,\gamma^\alpha_{qn}\,g^{mn}\,\gamma^p_{\alpha r}.
\endgathered\qquad
\mytag{4.19}
$$
Upon the formal change of summation indices $s\to n\to m\to h\to\alpha$
in the second term and $s\to r\to h$ in the third term respectively we
write \mythetag{4.19} as
$$
\gathered
\goth R^p_{qij}=\sum^3_{m=0}\sum^3_{n=0}\sum^3_{r=0}
\sum^4_{\alpha=1}\frac{L_{\boldsymbol\Upsilon_i}(\Gamma^r_{jm})
-L_{\boldsymbol\Upsilon_j}(\Gamma^r_{im})}{4}
\,\gamma^\alpha_{qn}\,g^{mn}\,\gamma^p_{\alpha r}-\\
\hskip 3em -\sum^3_{h=0}\sum^3_{m=0}\sum^3_{n=0}\sum^3_{r=0}
\sum^4_{\alpha=1}\frac{\Gamma^r_{jh}\,\Gamma^h_{i\kern 0.5pt m}
\,\gamma^\alpha_{qn}\,g^{mn}\,\gamma^p_{\alpha\kern 0.5pt r}}{4}\,+\\
\hskip -3em +\sum^3_{h=0}\sum^3_{m=0}\sum^3_{n=0}\sum^3_{r=0}
\sum^4_{\alpha=1}\frac{\Gamma^r_{i\kern 0.5pt h}\,\Gamma^h_{jm}
\,\gamma^\alpha_{qn}\,g^{mn}\,\gamma^p_{\alpha r}}{4}\,-\\
\hskip 3em -\sum^3_{k=0}\sum^3_{m=0}\sum^3_{n=0}\sum^3_{r=0}
\sum^4_{\alpha=1}\frac{c^{\,k}_{ij}}{4}\,\Gamma^r_{km}
\,\gamma^\alpha_{qn}\,g^{mn}\,\gamma^p_{\alpha r}.
\endgathered\qquad
\mytag{4.20}
$$
And finally, here is the last transformation that brings the formula 
\mythetag{4.20} to 
$$
\aligned
\goth R^p_{qij}&=\frac{1}{4}\sum^3_{m=0}\sum^3_{n=0}\sum^3_{r=0}
\sum^4_{\alpha=1}\left(L_{\boldsymbol\Upsilon_i}(\Gamma^r_{jm})
-L_{\boldsymbol\Upsilon_j}(\Gamma^r_{im})\,+
\vphantom{\sum^3_{h=0}}\right.\\
&\left.+\sum^3_{h=0}\left(\Gamma^r_{i\kern 0.5pt h}\,\Gamma^h_{jm}
-\Gamma^r_{jh}\,\Gamma^h_{i\kern 0.5pt m}\right)
-\sum^3_{k=0}\frac{c^{\,k}_{ij}}{4}\,\Gamma^r_{km}\right)
\gamma^\alpha_{qn}\ g^{mn}\,\gamma^p_{\alpha r}.
\endaligned\qquad
\mytag{4.21}
$$
Comparing \mythetag{4.21} with the formula \mythetag{4.7}, we can 
write the following ultimate result:
$$
\hskip -2em
\goth R^p_{qij}=\frac{1}{4}\sum^3_{m=0}\sum^3_{n=0}\sum^3_{r=0}
\sum^4_{\alpha=1}R^r_{mij}\,\gamma^\alpha_{qn}\ g^{mn}
\,\gamma^p_{\alpha r}.
\mytag{4.22}
$$
Note that the formula \mythetag{4.22} is quite similar to 
\mythetag{1.15}. However, like \mythetag{4.11} and \mythetag{4.12}
and unlike \mythetag{1.15}, it is a tensorial formula. For this reason,
being proved for some special pair of frames, it remains valid for 
an arbitrary pair of frames.
\mytheorem{4.1} The spinor curvature tensor $\eufb R$ of the 
torsion-free metric connection $(\Gamma,\Alpha,\bar{\Alpha})$ 
in the bundle of Dirac spinors $DM$ is related to the corresponding
Riemann curvature tensor $\bold R$ by means of the formula
\mythetag{4.22}. This formula is valid for any two frames 
$\boldsymbol\Psi_1,\,\boldsymbol\Psi_2,\,\boldsymbol\Psi_3,\,
\boldsymbol\Psi_4$ and $\boldsymbol\Upsilon_0,\,\boldsymbol
\Upsilon_1,\,\boldsymbol\Upsilon_2,\,\boldsymbol\Upsilon_3$
no matter holonomic or non-holonomic, special or not special,
and, if special, no matter being in frame association 
\mythetag{4.1} or not.
\endproclaim
\head
5. Conclusions.
\endhead
     The main result of this paper is that the formulas
\mythetag{1.14} and \mythetag{1.15} represent the spinor
components of the same metric connection in the bundle of
Dirac spinors. The formula \mythetag{1.14} is a general
formula, while \mythetag{1.15} is its specialization. The
formula \mythetag{1.15} is an important specialization
since, for instance, it is convenient for proving the 
formula \mythetag{4.22}.

\enddocument
\end